\documentclass[12pt]{article}
\usepackage{amsthm}
\usepackage{amsfonts}
\usepackage{amsmath}
\usepackage{amscd}
\usepackage{amssymb}
\usepackage{pstricks,pst-node,pst-tree}
\usepackage{stmaryrd}
\usepackage{newalg}
\theoremstyle{plain}
\newtheorem{thm}{Theorem}[section]
\newtheorem{lem}[thm]{Lemma}
\newtheorem{prop}[thm]{Proposition}

\newtheorem*{comp}{Computational Result}

\theoremstyle{definition}

\newtheorem{conj}{Conjecture}

\theoremstyle{remark}
\newtheorem*{rem}{Remark}

\renewcommand{\leq}{\leqslant}
\renewcommand{\geq}{\geqslant}
\title{REGULAR STEINHAUS GRAPHS\\ OF ODD DEGREE}
\date{02/09/2009}
\author{Jonathan \textsc{Chappelon}}
\begin{document}
\maketitle
\begin{abstract}
A Steinhaus matrix is a binary square matrix of size $n$ which is symmetric, with diagonal of zeros, and whose upper-triangular coefficients satisfy $a_{i,j}=a_{i-1,j-1}+a_{i-1,j}$ for all $2\leq i<j\leq n$. Steinhaus matrices are determined by their first row. A Steinhaus graph is a simple graph whose adjacency matrix is a Steinhaus matrix. We give a short new proof of a theorem, due to Dymacek, which states that even Steinhaus graphs, i.e. those with all vertex degrees even, have doubly-symmetric Steinhaus matrices. In 1979 Dymacek conjectured that the complete graph on two vertices $K_2$ is the only regular Steinhaus graph of odd degree. Using Dymacek's theorem, we prove that if $(a_{i,j})_{1\leq i,j\leq n}$ is a Steinhaus matrix associated with a regular Steinhaus graph of odd degree then its sub-matrix $(a_{i,j})_{2\leq i,j\leq n-1}$ is a multi-symmetric matrix, that is a doubly-symmetric matrix where each row of its upper-triangular part is a symmetric sequence. We prove that the multi-symmetric Steinhaus matrices of size $n$ whose Steinhaus graphs are regular modulo $4$, i.e. where all vertex degrees are equal modulo $4$, only depend on $\left\lceil \frac{n}{24}\right\rceil$ parameters for all even numbers $n$, and on $\left\lceil \frac{n}{30}\right\rceil$ parameters in the odd case. This result permits us to verify the Dymacek's conjecture up to $1500$ vertices in the odd case.
\end{abstract}
\section{Introduction}

Let $s=(a_1,a_2,\ldots,a_{n-1})$ be a binary sequence of length $n-1\geq1$ with entries $a_{j}$ in $\mathbb{F}_2=\{0,1\}$. The \textit{Steinhaus matrix} associated with $s$ is the square matrix $M(s)=(a_{i,j})$ of size $n$, defined as follows:
\begin{itemize}
\item
$a_{i,i}=0$\quad for all $1\leq i\leq n$,
\item
$a_{1,j}=a_{j-1}$\quad for all $2\leq j\leq n$,
\item
$a_{i,j}=a_{i-1,j-1}+a_{i-1,j}$\quad for all $2\leq i<j\leq n$,
\item
$a_{i,j}=a_{j,i}$\quad for all $1\leq i,j\leq n$.
\end{itemize}
By convention $M(\emptyset)=(0)$ is the Steinhaus matrix of size $n=1$ associated with the empty sequence. For example, the following matrix $M(s)$ in $\mathcal{M}_5(\mathbb{F}_2)$ is the Steinhaus matrix associated with the binary sequence $s=(1,1,0,0)$ of length $4$.
$$
M(s)=\left(
\begin{array}{ccccc}
0 & 1 & 1 & 0 & 0\\
1 & 0 & 0 & 1 & 0\\
1 & 0 & 0 & 1 & 1\\
0 & 1 & 1 & 0 & 0\\
0 & 0 & 1 & 0 & 0
\end{array}
\right)
$$
The set of all Steinhaus matrices of size $n\geq2$ will be denoted by $\mathcal{SM}_n(\mathbb{F}_2)$. It is clear that, for every positive integer $n$, the set $\mathcal{SM}_n(\mathbb{F}_2)$ has a cardinality of $2^{n-1}$.
\par The \textit{Steinhaus triangle} associated with $s$ is the upper-triangular part of the Steinhaus matrix $M(s)$. It was introduced by Hugo Steinhaus in 1963 \cite{Steinhaus1963}, who asked whether there exists a Steinhaus triangle containing as many $0$'s as $1$'s for each admissible size. Solutions of this problem appeared in \cite{Harborth1972,Eliahou2004}. A generalization of this problem to all finite cyclic groups was posed in \cite{Molluzzo1978} and was partially solved in \cite{Chappelon2008}.
\par The \textit{Steinhaus graph} associated with $s$ is the simple graph $G(s)$ on $n$ vertices whose adjacency matrix is the Steinhaus matrix $M(s)$. A vertex of a Steinhaus graph $G(s)$ is usually labelled by its corresponding row number in $M(s)$ and the $i$th vertex of $G(s)$ will be denoted by $V_i$. For instance, the following graph is the Steinhaus graph $G(s)$ associated with the sequence $s=(1,1,0,0)$.
$$
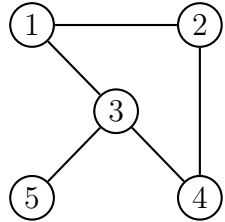
\psmatrix[colsep=0.5cm,rowsep=0.5cm,mnode=circle]
1&&2\\
&3\\
5&&4
\ncline{-}{1,1}{1,3}
\ncline{-}{1,1}{2,2}
\ncline{-}{1,3}{3,3}
\ncline{-}{3,3}{2,2}
\ncline{-}{3,1}{2,2}
\endpsmatrix
$$
For every positive integer $n$, the zero-edge graph on $n$ vertices is the Steinhaus graph associated with the sequence of zeros of length $n-1$.

Steinhaus graphs were introduced by Molluzzo in 1978 \cite{Molluzzo1978}. A general problem on Steinhaus graphs is that of characterizing those satisfying a given graph property. The bipartite Steinhaus graphs were characterized in \cite{Chang1999,Dymacek1986,Dymacek1995} and the planar ones in \cite{Dymacek2000a}. In \cite{Dymacek2000}, the following conjectures were made:

\begin{conj}\label{conj1}
The regular Steinhaus graphs of even degree are the zero-edge graph on $n$ vertices, for all positive integers $n$, and the Steinhaus graph $G(s)$ on $n=3m+1$ vertices generated by the periodic sequence $s=(1,1,0,\ldots,1,1,0)$ of length $3m$, for all positive integers $m$.
\end{conj}

\begin{conj}\label{conj2}
The complete graph on two vertices $K_2$ is the only regular Steinhaus graph of odd degree.
\end{conj}

These conjectures were verified up to $n\leq 25$ in 1988 by exhaustive search \cite{Bailey1988}. More recently \cite{Augier2008}, Augier and Eliahou extended the verification up to $n\leq117$ vertices by considering the weaker notion of \textit{parity-regular} Steinhaus graphs, i.e. Steinhaus graphs where all vertex degrees have the same parity. They searched regular graphs in the set of parity-regular Steinhaus graphs. This has enabled them to perform the verification because it is known that Steinhaus matrices associated with parity-regular Steinhaus graphs on $n$ vertices depend on approximately $n/3$ parameters \cite{Bailey1988,Augier2008}. This result is based on a theorem, due to Dymacek, which states that Steinhaus matrices associated with parity-regular Steinhaus graphs of even type are \textit{doubly-symmetric} matrices, i.e. where all the entries are symmetric with respect to the diagonal and the anti-diagonal of the matrices. A short new proof of this theorem is given in Section $2$. Using Dymacek's theorem, Bailey and Dymacek showed \cite{Bailey1988} that binary sequences associated with regular Steinhaus graphs of odd degree are of the form $(x_1,x_2,\ldots,x_k,x_k,\ldots,x_2,x_1,1)$. In Section $3$, we refine this result and, more precisely, we prove that if $(a_{i,j})_{1\leq i,j\leq n}$ is a Steinhaus matrix associated with a regular Steinhaus graph of odd degree, then its sub-matrix $(a_{i,j})_{2\leq i,j\leq n-1}$ is a \textit{multi-symmetric} Steinhaus matrix, i.e. a doubly-symmetric matrix where each row of the upper-triangular part is a symmetric sequence. A parametrization and a counting of multi-symmetric Steinhaus matrices of size $n$ are also given in Section $3$ for all $n\geq1$. In Section $4$, we show that, for a Steinhaus graph whose Steinhaus matrix is multi-symmetric, the knowledge of the vertex degrees modulo $4$ leads to a system of binary equations on the entries of its Steinhaus matrix. In Section $5$, we study the special case of multi-symmetric Steinhaus matrices whose Steinhaus graphs are \textit{regular modulo $4$}, i.e. where all vertex degrees are equal modulo $4$. We show that a such matrix of size $n$ only depends on $\left\lceil \frac{n}{24}\right\rceil$ parameters for all $n$ even, and on $\left\lceil \frac{n}{30}\right\rceil$ parameters in the odd case. Using these parametrizations, we obtain, by computer search, that for all positive integers $n\leq 1500$, the zero-edge graph on $n$ vertices is the only Steinhaus graph on $n$ vertices with a multi-symmetric matrix and which is regular modulo $4$. This permits us to extend the verification of Conjecture~\ref{conj2} up to $1500$ vertices.

\section{A new proof of Dymacek's theorem}

Recall that a square matrix $M=(a_{i,j})$ of size $n\geq1$ is said to be \textit{doubly-symmetric} if the entries of $M$ are symmetric with respect to the diagonal and to the anti-diagonal of $M$, that is
$$
a_{i,j} = a_{j,i} = a_{n-j+1,n-i+1},\quad \text{for\ all}\ 1\leq i,j\leq n.
$$

In \cite{Dymacek2000}, Dymacek characterized the parity-regular Steinhaus graphs. These results are based on the following theorem on parity-regular Steinhaus graphs of even type, where all vertex degrees are even.

\begin{thm}[Dymacek's theorem]\label{thm6}
The Steinhaus matrix of a parity-regular Steinhaus graph of even type is doubly-symmetric.
\end{thm}

In this section we give a new easier proof of Dymacek's theorem. The main idea of our proof is that the anti-diagonal entries of a Steinhaus matrix are determined by the vertex degrees of its associated Steinhaus graph.

\begin{thm}\label{thm1}
Let $G$ be a Steinhaus graph on $n\geq2$ vertices and $M=(a_{i,j})$ its associated Steinhaus matrix. Then every anti-diagonal entry of $M$ can be expressed by means of the vertex degrees of $G$. If we denote by $\deg(V_i)$ the degree of the vertex $V_i$ in $G$, then for all $1\leq i\leq\left\lfloor \frac{n}{2}\right\rfloor$, we have
$$
a_{i,n-i+1} \equiv \sum_{k=0}^{i-1}{\binom{i-1}{k}\deg\left(V_{i+k+1}\right)} \equiv \sum_{k=0}^{i-1}{\binom{i-1}{k}\deg\left(V_{n-i-k}\right)} \pmod{2}.
$$
\end{thm}

The proof is based on the following lemma which shows that each entry of the upper-triangular part of a Steinhaus matrix $M=(a_{i,j})$ can be expressed by means of the entries of the first row $\{a_{1,2} , \ldots , a_{1,n}\}$, the last column $\{a_{1,n} , \ldots , a_{n-1,n}\}$ or the over-diagonal $\{a_{1,2} , \ldots , a_{n-1,n}\}$ of $M$.

\begin{lem}\label{lem1}
Let $M=(a_{i,j})$ be a Steinhaus matrix of size $n\geq2$. Then, for all $1\leq i<j\leq n$, we have
$$
a_{i,j} = \sum_{k=0}^{i-1}{\binom{i-1}{k}a_{1,j-k}} = \sum_{k=0}^{n-j}{\binom{n-j}{k}a_{i+k,n}} = \sum_{k=0}^{j-i-1}{\binom{j-i-1}{k}a_{i+k,i+k+1}}.
$$
\end{lem}

\begin{proof}
Easily follows from the relation: $a_{i,j} = a_{i-1,j-1}+a_{i-1,j}$ for all $2\leq i<j\leq n$.
\end{proof}

\begin{proof}[Proof of Theorem~\ref{thm1}]
We begin by expressing each vertex degree of the Steinhaus graph $G$ by means of the entries of the first row, the last column and the over-diagonal of $M$. Here we view the entries $a_{i,j}$ as $0$,$1$ integers. For all $2\leq i\leq n-1$, we obtain
$$
\deg(V_i)
\begin{array}[t]{l}
= \displaystyle\sum_{j=1}^{n}{a_{i,j}} = \sum_{j=1}^{i-1}{a_{j,i}} + \sum_{j=i+1}^{n}{a_{i,j}}\\
\equiv \displaystyle\sum_{j=1}^{i-1}{(a_{j,i+1}+a_{j+1,i+1})} + \sum_{j=i+1}^{n}{(a_{i-1,j-1}+a_{i-1,j})}\\
\equiv \displaystyle\sum_{j=1}^{i-1}{a_{j,i+1}} + \sum_{j=2}^{i}{a_{j,i+1}} + \sum_{j=i}^{n-1}{a_{i-1,j}} + \sum_{j=i+1}^{n}{a_{i-1,j}}\\
\equiv a_{1,i+1} + a_{i,i+1} + a_{i-1,i} + a_{i-1,n} \pmod{2}.
\end{array}
$$
By Lemma~\ref{lem1}, it follows that
$$
\sum_{k=0}^{i-1}{\binom{i-1}{k}\deg\left(V_{i+k+1}\right)}
\begin{array}[t]{l}
\equiv \displaystyle\sum_{k=0}^{i-1}{\binom{i-1}{k}(a_{1,i+k+2}+a_{i+k+1,i+k+2}+a_{i+k,i+k+1}+a_{i+k,n})}\\
\equiv \displaystyle\sum_{k=0}^{i-1}{\binom{i-1}{k}a_{1,2i-k+1}} + \sum_{k=0}^{i-1}{\binom{i-1}{k}a_{i+k+1,i+k+2}}\\
+ \displaystyle\sum_{k=0}^{i-1}{\binom{i-1}{k}a_{i+k,i+k+1}} + \sum_{k=0}^{i-1}{\binom{i-1}{k}a_{i+k,n}}\\
\equiv a_{i,2i+1}+a_{i+1,2i+1}+a_{i,2i}+a_{i,n-i+1}\equiv a_{i,n-i+1} \pmod{2},
\end{array}
$$
for all $1\leq i\leq\left\lfloor \frac{n}{2}\right\rfloor$. The second congruence can be treated by the same way.
\end{proof}

\begin{rem}
We deduce from Theorem~\ref{thm1} a necessary condition on the vertex degrees of a given labelled graph to be a Steinhaus graph. Indeed, vertex degrees of a Steinhaus graph on $n$ vertices must satisfy the following binary equations:
$$
\sum_{k=0}^{i-1}{\binom{i-1}{k}\deg\left(V_{i+k+1}\right)} \equiv \sum_{k=0}^{i-1}{\binom{i-1}{k}\deg\left(V_{n-i-k}\right)} \pmod{2},\quad \text{for\ all}\ 1\leq i\leq\left\lfloor \frac{n}{2}\right\rfloor.
$$
More generally, an open problem, corresponding to Question $3$ in \cite{Dymacek1996}, is to determine if an arbitrary graph, not necessary labelled, is isomorphic to a Steinhaus graph.
\end{rem}

Now, we characterize doubly-symmetric Steinhaus matrices.

\begin{prop}\label{prop1}
Let $M=(a_{i,j})$ be a Steinhaus matrix of size $n\geq3$. Then the following assertions are equivalent:
\begin{itemize}
\item[(i)]
the matrix $M$ is doubly-symmetric,
\item[(ii)]
the over-diagonal of $M$ is a symmetric sequence,
\item[(iii)]
the entries $a_{i,n-i+1}$ of the anti-diagonal of $M$ vanish for all $1\leq i\leq \left\lfloor \frac{n-1}{2}\right\rfloor$.
\end{itemize}
\end{prop}

\begin{proof}
\ \\[2ex]
$\mathbf{(i)\Longrightarrow(ii):}$ Trivial.\\[2ex]
$\mathbf{(ii)\Longrightarrow(iii):}$ Suppose that the over-diagonal of $M$ is a symmetric sequence, that is
$$
a_{i,i+1} = a_{n-i,n-i+1},
$$
for all $1\leq i\leq n-1$. If $n$ is odd, then we have
$$
a_{i,n-i+1} = \displaystyle\sum_{k=0}^{n-2i}{\binom{n-2i}{k}a_{i+k,i+k+1}} = \displaystyle\sum_{k=0}^{\frac{n-2i+1}{2}}{\binom{n-2i}{k}(a_{i+k,i+k+1}+a_{n-i-k,n-i-k+1})} = 0,
$$
for all $1\leq i\leq\left\lfloor \frac{n-1}{2}\right\rfloor$. Otherwise, if $n$ is even, then we obtain
$$
a_{i,n-i+1} = \displaystyle\sum_{k=0}^{\frac{n}{2}-i-1}{\binom{n-2i}{k}(a_{i+k,i+k+1}+a_{n-i-k,n-i-k+1})} + 2\binom{n-2i-1}{\frac{n}{2}-i}a_{\frac{n}{2},\frac{n}{2}+1} = 0,
$$
for all $1\leq i\leq\left\lfloor \frac{n-1}{2}\right\rfloor$.\\[2ex]
$\mathbf{(iii)\Longrightarrow(i):}$ By induction on $n$. Consider the sub-matrix $N=(a_{i,j})_{2\leq i,j\leq n-1}$ that is a Steinhaus matrix of size $n-2$. By induction hypothesis, the matrix $N$ is doubly-symmetric. Then it remains to prove that $a_{1,j}=a_{n-j+1,n}$ for all $2\leq j\leq n$. First, since $a_{1,n}=0$, it follows that $a_{1,n-1}=a_{1,n}+a_{2,n}=a_{2,n}$ and for all $2\leq j\leq n-2$, we have
$$
a_{1,j} = \sum_{k=j+1}^{n-1}{a_{2,k}} + a_{1,n-1} = \sum_{k=2}^{n-j}{a_{k,n-1}} + a_{2,n} = a_{n-j+1,n}.
$$
\end{proof}

We are now ready to prove Dymacek's theorem.

\begin{proof}[Proof of Theorem~\ref{thm6}]
Let $G$ be a parity-regular Steinhaus graph of even type on $n$ vertices and $M=(a_{i,j})$ its Steinhaus matrix. If $n=1$, then $M=(0)$ which is trivially doubly-symmetric. Otherwise, for $n\geq2$, Theorem~\ref{thm1} implies that
$$
a_{i,n-i+1} \equiv \sum_{k=0}^{i-1}{\binom{i-1}{k}\deg\left(V_{i+k+1}\right)} \equiv 0 \pmod{2},
$$
for all $1\leq i\leq \left\lfloor \frac{n}{2}\right\rfloor$. Finally, the matrix $M$ is doubly-symmetric by Proposition~\ref{prop1}.
\end{proof}

\section{Multi-symmetric Steinhaus matrices}

In this section, we will study in detail the structure of Steinhaus matrices associated with regular Steinhaus graphs of odd degree.

Let $G$ be a Steinhaus graph on $n\geq1$ vertices. Then, for every integer $1\leq i\leq n$, we denote by $G\setminus\{V_i\}$ the graph obtained from $G$ by deleting its $i$th vertex $V_i$ and its incident edges in $G$. Since the adjacency matrix of the graph $G\setminus\{V_1\}$ (resp. $G\setminus\{V_n\}$) is the Steinhaus matrix obtained by removing the first row (resp. the last column) in the adjacency matrix of $G$, it follows that the graph $G\setminus\{V_1\}$ (resp. $G\setminus\{V_n\}$) is a Steinhaus graph on $n-1$ vertices.

Bailey and Dymacek studied the regular Steinhaus graphs of odd degree in \cite{Bailey1988}, where the following theorem is stated, using Dymacek's theorem.

\begin{thm}[\cite{Bailey1988}]\label{thm7}
Let $G$ be a regular Steinhaus graph of odd degree $d$ on $2n\geq4$ vertices. Then $d=n$, the Steinhaus graph $G\setminus\{V_1,V_{2n}\}$ is regular of even degree $n-1$, and $a_{1,j}=a_{1,2n-j+1}$ for all $2\leq j\leq 2n-1$.
\end{thm}

\begin{rem}
In every simple graph, there are an even number of vertices of odd degree. Therefore parity-regular Steinhaus graphs of odd type and thus regular Steinhaus graphs of odd degree have an even number of vertices.
\end{rem}

In their theorem, the authors studied the form of the sequence associated with $G$. We are more interested in the Steinhaus matrix of $G\setminus\{V_1,V_{2n}\}$ in the sequel.

Recall that a square matrix of size $n\geq1$ is said to be \textit{multi-symmetric} if $M$ is doubly-symmetric and each row of the upper-triangular part of $M$ is a symmetric sequence, that is
$$
a_{i,j} = a_{i,n-j+i+1},\quad \text{for\ all}\ 1\leq i<j\leq n.
$$

First, it is easy to see that each column of the upper-triangular part of a multi-symmetric matrix is also a symmetric sequence.

\begin{prop}\label{prop3}
Let $M=(a_{i,j})$ be a multi-symmetric matrix of size $n$. Then, each column of the upper-triangular part of $M$ is a symmetric sequence, that is $a_{i,j} = a_{j-i,j}$ for all $1\leq i<j\leq n$.
\end{prop}

\begin{proof}
Easily follows from the relation: $a_{i,j} = a_{i,n-j+i+1} = a_{j-i,n-i+1} = a_{j-i,j}$ for all $1\leq i<j\leq n$.
\end{proof}

As for doubly-symmetric Steinhaus matrices, multi-symmetric Steinhaus matrices can be characterized as follows.

\begin{prop}\label{prop4}
Let $M=(a_{i,j})$ be a Steinhaus matrix of size $n\geq3$. Then the following assertions are equivalent:
\begin{itemize}
\item[(i)]
the matrix $M$ is multi-symmetric,
\item[(ii)]
the first row, the last column and the over-diagonal of $M$ are symmetric sequences,
\item[(iii)]
the entries $a_{i,n-i+1}$, $a_{n-2i+1,n-i+1}$ and $a_{i,2i}$ vanish for all $1\leq i\leq\left\lfloor\frac{n-1}{2} \right\rfloor$.
\end{itemize}
\end{prop}

\begin{proof}
Similar to the proof of Proposition~\ref{prop1} and by using Lemma~\ref{lem1} and Proposition~\ref{prop3}.
\end{proof}

We now refine Theorem~\ref{thm7}.

\begin{thm}\label{thm5}
Let $G$ be a regular Steinhaus graph of odd degree $n$ on $2n\geq4$ vertices. Then $G\setminus\{V_1,V_{2n}\}$ is a regular Steinhaus graph of even degree $n-1$ whose associated Steinhaus matrix is multi-symmetric.
\end{thm}

\begin{proof}
Let $M=(a_{i,j})$ be the Steinhaus matrix associated with $G$. Theorem~\ref{thm7} implies that the Steinhaus graph $G\setminus\{V_1,V_{2n}\}$ is regular of even degree $n-1$ and that we have
$$
a_{1,j}=a_{1,2n-j+1},
$$
for all $2\leq j\leq 2n-1$. Therefore, for all $3\leq j\leq 2n-1$, we have
$$
a_{2,j} + a_{2,2n-j+2} = (a_{1,j-1} + a_{1,j}) + (a_{1,2n-j+1} + a_{1,2n-j+2}) = (a_{1,j-1} + a_{1,2n-j+2}) + (a_{1,j} + a_{1,2n-j+1}) = 0.
$$
Then the first row of the matrix $B=(a_{i,j})_{2\leq i,j\leq 2n-1}$, the Steinhaus matrix of the graph $G\setminus\{V_1,V_{2n}\}$, is a symmetric sequence. Moreover, by Dymacek's theorem, the matrix $B$ is doubly-symmetric. Finally, by Proposition~\ref{prop4}, the matrix $B$ is multi-symmetric.
\end{proof}

\begin{rem}
By Theorem~\ref{thm5}, it is easy to show that Conjecture~\ref{conj1} implies Conjecture~\ref{conj2}. Indeed, if Conjecture~\ref{conj1} is true, then the zero-edge graph on $n$ vertices is the only regular Steinhaus graph of even degree whose Steinhaus matrix is multi-symmetric. It follows, by Theorem~\ref{thm5}, that if $G(s)$ is a regular Steinhaus graph of odd degree on $n+2$ vertices then $s=(0,\ldots,0,1)$ or $s=(1,\ldots,1)$. Therefore the Steinhaus graph $G(s)$ is the star graph on $n+2$ vertices which is not a regular Steinhaus graph.
\end{rem}

In the sequel of this section we will study in detail the multi-symmetric Steinhaus matrices. First, in order to determine a parametrization of these matrices, we introduce the following operator
$$
T : \mathcal{SM}_n(\mathbb{F}_2) \longrightarrow \mathcal{SM}_{n-3}(\mathbb{F}_2),
$$
which assigns to each matrix $M=(a_{i,j})$ in $\mathcal{SM}_n(\mathbb{F}_2)$ the Steinhaus matrix $T(M)=(b_{i,j})$ in $\mathcal{SM}_{n-3}(\mathbb{F}_2)$ defined by $b_{i,j} = a_{i-1,j-2}$, for all $1\leq i<j\leq n-3$. As depicted in the following matrix, the upper-triangular part of $M$ is an extension of the upper-triangular part of $T(M)$.
$$
\left(
\begin{array}{ccccccccccccc}
0 & a_{1,2} & a_{1,3} & a_{1,4} & a_{1,5} & a_{1,6} & \cdots & \cdots & a_{1,n-4} & a_{1,n-3} & a_{1,n-2} & a_{1,n-1} & a_{1,n}\\
& 0 & a_{2,3} & \mathbf{b_{1,2}} & \mathbf{b_{1,3}} & \mathbf{b_{1,4}} & \cdots & \cdots & \cdots & \mathbf{b_{1,n-5}} & \mathbf{b_{1,n-4}} & \mathbf{b_{1,n-3}} & a_{2,n}\\
& & 0 & a_{3,4} & \mathbf{b_{2,3}} & \mathbf{b_{2,4}} & & & & & \mathbf{b_{2,n-4}} & \mathbf{b_{2,n-3}} & a_{3,n}\\
& & & 0 & a_{4,5} & \mathbf{b_{3,4}} & & & & & & \mathbf{b_{3,n-3}} & a_{4,n}\\
& & & & 0 & a_{5,6} & \ddots & & & & & \vdots & a_{5,n}\\
& & & & & 0 & \ddots & \ddots & & & & \vdots & \vdots\\
& & & & & & \ddots & \ddots & \ddots & & & \vdots & \vdots\\
& & & & & & & 0 & a_{n-5,n-4} & \mathbf{b_{n-6,n-5}} & \mathbf{b_{n-6,n-4}} & \mathbf{b_{n-6,n-3}} & a_{n-5,n}\\
& & & & & & & & 0 & a_{n-4,n-3} & \mathbf{b_{n-5,n-4}} & \mathbf{b_{n-5,n-3}} & a_{n-4,n}\\
& & & & & & & & & 0 & a_{n-3,n-2} & \mathbf{b_{n-4,n-3}} & a_{n-3,n}\\
& & & & & & & & & & 0 & a_{n-2,n-1} & a_{n-2,n}\\
& & & & & & & & & & & 0 & a_{n-1,n}\\
& & & & & & & & & & & & 0
\end{array}
\right)
$$

\begin{prop}\label{prop11}
Let $M=(a_{i,j})$ be a Steinhaus matrix of size $n\geq4$. Then the extension $M$ of $T(M)$ only depends on the parameters $a_{1,2}$, $a_{1,j_0}$ and $a_{1,n}$, with $j_0$ in $\{3,\ldots,n-1\}$.
\end{prop}

\begin{proof}
Let $3\leq j_0\leq n-1$. Each entry $a_{1,j}$, for $3\leq j\leq n-1$, can be expressed by means of $a_{1,j_0}$ and the entries of $T(M)=(b_{i,j})$. Indeed, we have
$$
\begin{array}{rll}
a_{1,j} & = a_{1,j_0} + \displaystyle\sum_{k=j-1}^{j_0-2}{b_{1,k}}, & \text{for\ all}\ 3\leq j<j_0,\\
a_{1,j} & = a_{1,j_0} + \displaystyle\sum_{k=j_0-1}^{j-2}{b_{1,k}}, & \text{for\ all}\ j_0 < j\leq n-1.
\end{array}
$$
Then the entries $a_{1,2}$, $a_{1,j_0}$ and $a_{1,n}$ determine the extension $M$ of $T(M)$.
\end{proof}

Therefore, for every Steinhaus matrix $N$ of size $n-3$, there exist $8$ distinct Steinhaus matrices $M$ of size $n$ such that $T(M)=N$. We can also use this operator to determine parametrizations of multi-symmetric Steinhaus matrices.

\begin{prop}\label{prop7}
Let $M=(a_{i,j})$ be a multi-symmetric Steinhaus matrix of size $n$. Let $j_i$ be an element of the set $\{2i+1,\ldots,n-i\}$ for all $1\leq i\leq\left\lfloor \frac{n-1}{3}\right\rfloor$. Then the matrix $M$ depends on the following parameters:
\begin{itemize}
\item $a_{1,j_1}$ and $\left\{a_{2i,j_{2i}}\ \middle|\ 1\leq i\leq\left\lceil \frac{n}{6}\right\rceil-1\right\}$, for $n$ even,
\item $\left\{a_{2i+1,j_{2i+1}}\ \middle|\ 0\leq i\leq\left\lceil \frac{n-3}{6}\right\rceil-1 \right\}$, for $n$ odd.
\end{itemize}
\end{prop}

\begin{proof}
Let $M=(a_{i,j})$ be a multi-symmetric matrix of size $n$. We consider the sub-matrices $T(M),\ T^2(M)=T(T(M)),\ T^3(M),\ T^4(M),\ \ldots$. By successive application of Proposition~\ref{prop11} on the extension $T^{i-1}(M)$ of $T^i(M)$ and since the entries $a_{i,n-i+1}$, $a_{n-2i+1,n-i+1}$ and $a_{i,2i}$ vanish for all $1\leq i\leq\left\lfloor \frac{n-1}{2}\right\rfloor$ by Proposition~\ref{prop4}, the parametrizations of the multi-symmetric matrix $M$ follow.
\end{proof}

For all positive integers $n$, the number of multi-symmetric Steinhaus matrices of size $n$ immediately follows.

\begin{thm}
Let $n$ be a positive integer. If we denote by $MS(n)$ the number of multi-symmetric Steinhaus matrices of size $n$, then we have
$$
MS(n) = \left\{
\begin{array}{ll}
2^{\left\lceil \frac{n}{6}\right\rceil} & ,\ \text{for}\ n\ \text{even},\\
2^{\left\lceil \frac{n-3}{6}\right\rceil} & ,\ \text{for}\ n\ \text{odd}.
\end{array}
\right.
$$
\end{thm}

\section{Vertex degrees of Steinhaus graphs associated with\\ multi-symmetric Steinhaus matrices}

In this section, we analyse the vertex degrees of a Steinhaus graph associated with a multi-symmetric Steinhaus matrix of size $n$. We begin with the case of doubly-symmetric Steinhaus matrices.

\begin{prop}\label{prop10}
Let $n$ be a positive integer and $G$ be a Steinhaus graph on $n$ vertices whose Steinhaus matrix is doubly-symmetric. Then, for all $1\leq i\leq n$, we have
$$
\deg(V_i) = \deg(V_{n-i+1}).
$$
\end{prop}

\begin{proof}
If we denote by $M=(a_{i,j})$ the Steinhaus matrix associated with the graph $G$, then, for all $1\leq i\leq n$, we have
$$
\deg(V_i) = \sum_{j=1}^{n}{a_{i,j}} = \sum_{j=1}^{n}{a_{n-j+1,n-i+1}} = \sum_{j=1}^{n}{a_{j,n-i+1}} = \deg(V_{n-i+1}).
$$
\end{proof}

We shall now see that, for a Steinhaus graph associated with a multi-symmetric Steinhaus matrix, the knowledge of the vertex degrees modulo $4$ imposes strong conditions on the entries of its Steinhaus matrix. In order to prove this result, we distinguish different cases depending on the parity of $n$.

\begin{prop}\label{prop6}
Let $n$ be an even number and $G$ be a Steinhaus graph on $n$ vertices whose Steinhaus matrix $M=(a_{i,j})$ is multi-symmetric. Then, we have
$$
\begin{array}{l}
\deg(V_1) = \deg(V_n) \equiv a_{1,\frac{n}{2}+1} \pmod{2},\\
\deg(V_2) = \deg(V_{n-1}) \equiv 2a_{1,\frac{n}{2}+1} \pmod{4},\\
\deg(V_3) = \deg(V_{n-2}) \equiv 2a_{2,\frac{n}{2}+1} \pmod{4},\\
\deg(V_{2i}) = \deg(V_{n-2i+1}) \equiv 2a_{2,2i+1} + 2a_{i,2i+1} \pmod{4},\quad \text{for\ all}\ 2\leq i\leq\frac{n}{2}-2.
\end{array}
$$
\end{prop}

\begin{proof}
First, Proposition~\ref{prop4} implies that the entries $a_{i,2i}$ and $a_{2i+1,\frac{n}{2}+i+1}$ vanish for all $1\leq i\leq\frac{n}{2}-1$. This leads to
$$
\begin{array}{l}
\deg(V_1) = \displaystyle\sum_{j=2}^{n}{a_{1,j}} = \sum_{j=2}^{\frac{n}{2}}{(a_{1,j}+a_{1,n-j+2})}+a_{1,\frac{n}{2}+1} \equiv a_{1,\frac{n}{2}+1} \pmod{2},\\
\deg(V_2) = a_{1,2} + \displaystyle\sum_{j=3}^{\frac{n}{2}+1}{(a_{2,j}+a_{2,n-j+3})} = 2\sum_{j=3}^{\frac{n}{2}+1}{a_{2,j}} \equiv 2a_{1,2} + 2a_{1,\frac{n}{2}+1} \equiv 2a_{1,\frac{n}{2}+1} \pmod{4},\\
\deg(V_3) = (a_{1,3}+a_{2,3})+\displaystyle\sum_{j=4}^{\frac{n}{2}+1}{(a_{3,j}+a_{3,n-j+4})}+a_{3,\frac{n}{2}+2} = 2a_{2,3}+2\sum_{j=4}^{\frac{n}{2}+1}{a_{3,j}} \equiv 2a_{2,\frac{n}{2}+1} \pmod{4},
\end{array}
$$
and, for all $2\leq i\leq\frac{n}{2}-2$, we have
$$
\deg(V_{2i})
\begin{array}[t]{l}
= \displaystyle\sum_{j=i+1}^{2i-1}{(a_{j,2i}+a_{2i-j,2i})}+a_{i,2i}+\sum_{j=2i+1}^{\frac{n}{2}+i}{(a_{2i,j}+a_{2i,n-j+2i+1})}\\
= 2\displaystyle\sum_{j=i+1}^{2i-1}{a_{j,2i}} + 2\sum_{j=2i+1}^{\frac{n}{2}+i}{a_{2i,j}}\\
\equiv 2\displaystyle\sum_{j=i+1}^{2i-1}{a_{j,2i+1}} + 2\sum_{j=i+2}^{2i}{a_{j,2i+1}} + 2\sum_{j=2i}^{\frac{n}{2}+i-1}{a_{2i-1,j}} + 2\sum_{j=2i+1}^{\frac{n}{2}+i}{a_{2i-1,j}}\\
\equiv 2a_{i+1,2i+1} + 2a_{2i,2i+1} + 2a_{2i-1,2i} + 2a_{2i-1,\frac{n}{2}+i}\\
\equiv 2a_{i+1,2i+1} + 2a_{2i-1,2i+1} \equiv 2a_{2,2i+1} + 2a_{i,2i+1} \pmod{4}.
\end{array}
$$
Finally, we complete the proof by Proposition~\ref{prop10}.
\end{proof}

\begin{rem}
Let $n$ be an even number. In every Steinhaus graph on $n$ vertices whose Steinhaus matrix is multi-symmetric the fourth vertex $V_4$ has a degree divisible by $4$.
\end{rem}

\begin{prop}\label{prop8}
Let $n$ be an odd number and $G$ be a Steinhaus graph on $n$ vertices whose Steinhaus matrix $M=(a_{i,j})$ is multi-symmetric. Then, we have
$$
\begin{array}{l}
\deg(V_1) = \deg(V_n) \equiv 0 \pmod{2},\\
\deg(V_2) = \deg(V_{n-1}) \equiv 2a_{1,\frac{n+1}{2}} \pmod{4},\\
\deg(V_{2i}) \equiv 2a_{i+1,2i+1} + 2a_{2i-1,2i+1} + 2a_{2i-1,\frac{n-1}{2}+i} \pmod{4},\quad \text{for\ all}\ 2\leq i\leq\frac{n-3}{2},\\
\deg(V_{2i+1}) \equiv 2a_{2,2i+2} \pmod{4},\quad \text{for\ all}\ 1\leq i\leq\frac{n-3}{2}.
\end{array}
$$
\end{prop}

\begin{proof}
Proposition~\ref{prop4} implies that the entries $a_{i,2i}$ and $a_{2i,(n+1)/2+i}$ vanish for all $1\leq i \leq \frac{n-3}{2}$. Since each row and each column of the upper triangular part of $M$ is symmetric, we can use the relation
$$
\sum_{k=1}^{m}a_{i,j+k} \equiv a_{i-1,j} + a_{i-1,j+m} \pmod{2},\quad \text{for\ all}\ 2\leq i<j\leq n-m+1
$$
as in the proof of Proposition~\ref{prop6}, and the results follow.
\end{proof}

\begin{rem}
Let $n$ be an odd number. In every Steinhaus graph on $n$ vertices whose Steinhaus matrix is multi-symmetric the third vertex $V_3$ has a degree divisible by $4$.
\end{rem}

\section{Multi-symmetric Steinhaus matrices\\ of Steinhaus graphs with regularity modulo $4$}

In this section, we consider the multi-symmetric Steinhaus matrices associated with Steinhaus graphs which are regular modulo $4$, i.e. where all vertex degrees are equal modulo $4$. First, we determine an upper bound of the number of these matrices. Two cases are distinguished, according to the parity of $n$.

\begin{thm}\label{thm4}
For all odd numbers $n$, there are at most $2^{\left\lceil \frac{n}{30}\right\rceil}$ multi-symmetric Steinhaus matrices of size $n$ whose associated Steinhaus graphs are regular modulo $4$.
\end{thm}

\begin{proof}
Let $n$ be an odd number and $M=(a_{i,j})$ a multi-symmetric Steinhaus matrix of size $n$. By Proposition~\ref{prop7}, the matrix $M$ depends on the parameters $a_{2i+1,\frac{n+1}{2}+i}$ for $0\leq i\leq \left\lceil \frac{n-3}{6}\right\rceil-1$. If the Steinhaus graph associated with $M$ is regular modulo $4$, then Proposition~\ref{prop8} implies that $a_{2,2j}=0$ for all $2\leq j\leq \frac{n-1}{2}$ and thus
$$
a_{2i,2j} = \sum_{k=0}^{i-1}a_{2,2j-2k} =0,
$$
for all $1\leq i<j\leq \frac{n-1}{2}$.
\par If $n\equiv 1\pmod{4}$, then $\frac{n+1}{2}$ is odd and
$$
a_{4i+1,\frac{n+1}{2}+2i} = a_{4i,\frac{n-1}{2}+2i} + a_{4i,\frac{n+1}{2}+2i} = 0,
$$
for all $0\leq i\leq \left\lfloor \frac{\left\lceil \frac{n-3}{6}\right\rceil-1}{2}\right\rfloor$. Therefore the matrix $M$ can be parametrized by
$$
\left\{ a_{4i+3,\frac{n+3}{2}+2i}\ \middle|\ 0\leq i\leq m-1\right\},
$$
with
$$
m=\left\lceil \frac{\left\lceil \frac{n-3}{6}\right\rceil-1}{2}\right\rceil.
$$
Suppose that we know the $p$ parameters in
$$
P = \left\{ a_{4i+3,\frac{n+3}{2}+2i}\ \middle|\ m-p\leq i\leq m-1\right\}.
$$
Then, by Proposition~\ref{prop7} again, the multi-symmetric matrix $T^{4(m-p)-1}(M)$ can be parametrized by $P$. Therefore the entries
$$
\left\{ a_{i,2i+1}\ \middle|\ 4(m-p)\leq i\leq \frac{n-1}{2}-2(m-p)\right\}
$$
in $T^{4(m-p)-1}(M)$ depend on the parameters in $P$. Moreover, if the Steinhaus graph associated with $M$ is regular modulo $4$, then Proposition~\ref{prop8} implies that
$$
a_{2,2i+1} = a_{2i-1,2i+1} \equiv a_{i+1,2i+1} + a_{2i-1,\frac{n-1}{2}+i} \equiv a_{i+1,2i+1} + a_{\left(\frac{n+1}{2}-i\right)+1,2\left(\frac{n+1}{2}-i\right)+1} \pmod{2},
$$
for all $1\leq i\leq\frac{n-1}{2}$. If the inequality
$$
\frac{n+1}{2}-4(m-p)\geq4(m-p)
$$
holds, then the entries $a_{2,2i+1}$ depend on the parameters in $P$ for all $4(m-p)\leq i\leq \frac{n+1}{2}-4(m-p)$. Since we have $a_{2,2i}=0$ for all $4(m-p)\leq i\leq \frac{n+3}{2}-4(m-p)$, it follows that the entries
$$
\left\{ a_{i,j}\ \middle|\ 
\begin{array}{c}
2 \leq i \leq n+5-16(m-p)\\
8(m-p)+i-1 \leq j \leq n+3-8(m-p)
\end{array} \right\}
$$
depend on the parameters in $P$. Suppose now that $p$ is solution of the following inequality
$$
n+5-16(m-p) \geq 4(m-p)-1.
$$
Therefore the extension $M$ of $T^{(4(m-p)-1)}(M)$ depends on the entries $a_{i,n+3-8(m-p)}$ for $2\leq i\leq 4(m-p)-1$ and $a_{1,\frac{n+1}{2}}$ which vanishes by Proposition~\ref{prop8}. Thus, all the entries of the matrix $M$ depend on the $p$ parameters in $P$. Finally, a solution of this inequality can be obtained when
$$
p=\left\lceil \frac{n}{30}\right\rceil \geq \left\lceil \frac{\left\lceil \frac{n-3}{6}\right\rceil-1}{2}\right\rceil-\frac{n+6}{20}.
$$
\par If $n\equiv3\pmod{4}$, then $\frac{n+1}{2}$ is even and
$$
a_{4i+3,\frac{n+3}{2}+2i} = a_{4i+2,\frac{n+1}{2}+2i} + a_{4i+2,\frac{n+3}{2}+2i} = 0,
$$
for all $0\leq i\leq \left\lceil \frac{\left\lceil \frac{n-3}{6}\right\rceil-1}{2}\right\rceil-1$. Therefore the matrix $M$ can be parametrized by
$$
\left\{ a_{4i+1,\frac{n+1}{2}+2i}\ \middle|\ 0\leq i\leq m\right\}
$$
with
$$
m=\left\lfloor \frac{\left\lceil \frac{n-3}{6}\right\rceil-1}{2}\right\rfloor.
$$
As above, in the case $n\equiv1\pmod{4}$, we can prove that all the entries of the matrix $M$ depend on the $p$ parameters in
$$
\left\{ a_{4i+1,\frac{n+1}{2}+2i}\ \middle|\ m-p+1\leq i\leq m\right\}
$$
if $p$ is solution of the following inequality
$$
n-16(m-p)-4 \geq 4(m-p)+1.
$$
A solution is obtained when
$$
p = \left\lceil \frac{n}{30}\right\rceil \geq \left\lfloor \frac{\left\lceil \frac{n-3}{6}\right\rceil-1}{2}\right\rfloor - \frac{n-5}{20}.
$$
\end{proof}

\begin{thm}
For all even numbers $n$, there are at most $2^{\left\lceil \frac{n}{24}\right\rceil}$ multi-symmetric Steinhaus matrices of size $n$ whose associated Steinhaus graphs are regular modulo $4$.
\end{thm}

\begin{proof}[Sketch of proof]
Similar to the proof of Theorem~\ref{thm4}. Let $M=(a_{i,j})$ be a multi-symmetric Steinhaus matrix of even size $n$. First, by Proposition~\ref{prop7}, for all positive integers $p<m-1$ with $m=\left\lceil \frac{n}{6}\right\rceil$, the multi-symmetric Steinhaus matrix $T^{2\left(m-p-1\right)}(M)$ can be parametrized by the $p$ entries in
$$
P = \left\{ a_{2i,4i+1}\ \middle|\ m - p \leq i\leq m-1\right\}.
$$
Moreover, if the Steinhaus graph associated with $M$ is regular modulo $4$, then Proposition~\ref{prop6} implies that $a_{1,\frac{n}{2}+1}=0$ and $a_{2,2i+1} = a_{i,2i+1}$ for all $2\leq i\leq \frac{n}{2}-1$. It follows that the entries $a_{2i,n-2(m-p)+1}$ also depends on the parameters in $P$ for all $1\leq i\leq \frac{n}{2}-3(m-p)+2$. Finally, we can see that, if $p$ is solution of the following inequality
$$
\frac{n}{2}-3(m-p)+2 \geq m-p-1,
$$
then, as in the proof of Proposition~\ref{prop7}, the extension $M$ of $T^{2(m-p-1)}(M)$ depends on the entries $a_{2i,n-2(m-p)+1}$ for $1\leq i\leq m-p-1$ and thus all the entries of the matrix $M$ can be expressed by means of the $p$ parameters in $P$. We conclude the proof by observing that the inequality is obtained when
$$
p = \left\lceil \frac{n}{24}\right\rceil \geq \left\lceil \frac{n}{6}\right\rceil - \frac{n+6}{8}.
$$
\end{proof}

Using these explicit parametrizations of the multi-symmetric Steinhaus matrices whose Steinhaus graphs are regular modulo $4$, we obtain the following result by computer search:

\begin{comp}
For all positive integers $n\leq 1500$, the zero-edge graph on $n$ vertices is the only Steinhaus graph on $n$ vertices with a multi-symmetric Steinhaus matrix and which is regular modulo $4$.
\end{comp}

This result can be easily proved for all odd numbers in the special case of regular Steinhaus graphs on $n$ vertices whose Steinhaus matrices are multi-symmetric.

\begin{thm}
For all odd numbers $n$, there is no regular Steinhaus graph on $n$ vertices whose Steinhaus matrix is multi-symmetric, except the zero-edge graph on $n$ vertices.
\end{thm}

\begin{proof}
Let $n$ be an odd number. Let $G$ be a regular Steinhaus graph on $n$ vertices and $M=(a_{i,j})$ its Steinhaus matrix. Then Proposition~\ref{prop8} implies that
$$
\deg(V_i) \equiv 0 \pmod{4},
$$
for all $1\leq i\leq n$ and
$$
a_{2,2i+2} = 0,
$$
for all $1\leq i\leq\frac{n-3}{2}$. If we denote by $\oplus$ the addition in $\mathbb{F}_2$ and $+$ the addition in the integers, then we obtain
$$
\deg(V_3)
\begin{array}[t]{l}
= a_{1,3} + a_{2,3} + \displaystyle\sum_{j=4}^{n}{a_{3,j}} = (a_{1,2}\oplus a_{2,3}+a_{2,3}) + \sum_{j=2}^{\frac{n-3}{2}}{(a_{3,2j+1}+a_{3,2j+2})} + 2a_{3,n}\\
= 2a_{2,3} + \displaystyle\sum_{j=2}^{\frac{n-3}{2}}{(a_{2,2j}\oplus a_{2,2j+1} + a_{2,2j+1}\oplus a_{2,2j+2})} + 2(a_{2,n-1}\oplus a_{2,n})\\
= 2\displaystyle\sum_{j=1}^{\frac{n-1}{2}}{a_{2,2j+1}} = 2(a_{1,2}+\sum_{j=3}^{n}{a_{2,j}}) = 2\times\deg(V_2).
\end{array}
$$
This leads to $\deg(V_i)=0$ for all $1\leq i\leq n$ and thus $G$ is the zero-edge graph on $n$ vertices.
\end{proof}

Finally, the above computational result permits us to extend the verification of Conjecture~\ref{conj2} up to $n\leq 1500$ vertices. Indeed, as proved in the remark following Theorem~\ref{thm5}, for a Steinhaus graph $G$ on $2n$ vertices, if $G\setminus\{V_1,V_{2n}\}$ is the zero-edge graph on $2n-2$ vertices, then $G$ is the star graph on $2n$ vertices which is not a regular graph. Therefore, by Theorem~\ref{thm5}, we obtain

\begin{thm}
There is no regular Steinhaus graph of odd degree on $2<n\leq 1500$ vertices.
\end{thm}

\section*{Acknowledgments}

The author would like to thank Shalom Eliahou for introducing him to the subject and for his help in preparing this paper.

\nocite{*}
\bibliographystyle{plain}
\bibliography{biblio}

\noindent Jonathan Chappelon\\
Laboratoire de Math\'ematiques Pures et Appliqu\'ees Joseph Liouville, FR CNRS 2956\\
Universit\'e du Littoral C\^ote d'Opale\\
50 rue F. Buisson, B.P. 699, F-62228 Calais Cedex, France\\
e-mail: jonathan.chappelon@lmpa.univ-littoral.fr

\end{document}